# COLLECTIVE DESIGN OF AN ONLINE MATH TEXTBOOK: WHEN INDIVIDUAL AND COLLECTIVE DOCUMENTATION WORKS MEET

Hussein Sabra / Luc Trouche

INRP & LIRDHIST (Université Lyon 1, France)

*This paper focuses on the documentation work in the mathematics teaching. We show up the individual and collective components of this documentation. We present a new theoretical framework, the documentational approach which seems adapted for studying this issue. We applied it on a particular project of the French teachers association Sésamath: the collective design of an online mathematics textbook. We present a methodology for observing both this collective project and the case of Anaïs, a teacher involved in this project. We study her professional interest in Sésamath, her didactical interactions with the other project members, and their effects on her professional knowledge.*

## 1. INTRODUCTION

*Digitalization* deeply changes conditions of thinking and sharing knowledge at each level of the society (Pédauque 2006). Its visible manifestations, in mathematics teaching (Hoyles & Lagrange 2009), are both the *profusion* of resources available via Internet and the *diversification* of technologies that could be used by the teachers (USB keys, Interactive whiteboard, calculators, software). These evolutions dramatically modify the conditions of professional practice (for *preparing* as well as for *doing* the teaching). They prompt new collective forms of teachers work: on the one hand digitalization makes *potentially* easier the sharing and exchange of resources between teachers, on the other hand, its complexity (dispersion of resources, rapid technological evolutions) makes *necessary* for teachers to help each other. These new forms of collective work among mathematics teachers change the conditions for their professional growth: *teams, networks and communities* appear as new opportunities for teachers to learn (Krainer & Woods 2008).

The most significant phenomenon of this trend, in France, is, for us, the appearing of Sésamath[1], a French online association of mathematics teachers, aiming to provide mathematics teachers with free online resources. For achieving this goal, Sésamath develops collaborative work of teachers, around common projects (Sabra 2009). We should want to address, in this article, two questions: why do teachers engage in Sésamath? How do they articulate their work on resources for individual purpose and for Sésamath purpose?

In the following section, we introduce a new theoretical approach, which seems necessary for addressing these issues. In the third section, we present our experimental field and specify our questions. We set out our methodology in the fourth section. In the fifth section we display some preliminary results. In conclusion, we propose some new questions that this research raises.

## 2. THEORETICAL FRAMEWORK

Our research lies on a new approach to teachers' work with *resources* and professional development, already presented in CERME 6 (Gueudet & Trouche 2009) and developed further (Gueudet, Pepin & Trouche, to appear): *the documentational approach of didactics*.

It is built on a distinction between *resources* and *documents*, extending the one, introduced by the *instrumental approach* (Rabardel 1995), between *artefact* and *instrument*. The choice of the word "resource", instead of artefact, aims at catching a great variety of things intervening in teachers work: a textbook, a piece of software, a student's worksheet, an Internet resource, a discussion with a colleague, etc. We call *documentation work* what a teacher needs to do for designing her teaching: looking for resources, integrating them in her personal *resource system*, implementing it in practice, sharing it with colleagues, renewing it taking into account various feedback etc.

A teacher draws on resource sets for her documentation work. A process of genesis (Fig. 1) takes place, producing what we call a *document*[2], made of resources and a scheme (i.e. an invariant organization of the activity to perform a type of task - here a task for preparing and performing a given teaching). Each scheme encapsulated professional knowledge, both shaping teacher's activity and permanently reshaped by this activity. Shulman (1986) proposed a categorization of teacher's professional knowledge. We are particularly interested in one of these categories that Grossman (1990) developed further: the Pedagogical Content Knowledge (PCK), defined as knowledge that a teacher develops to help her students in their learning.

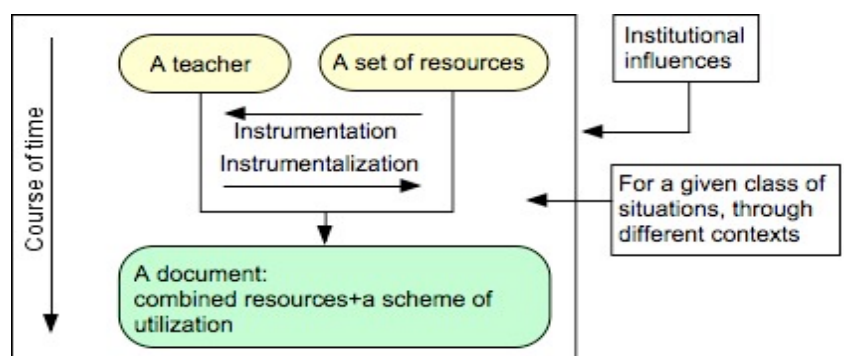

**Figure 1.** *Schematic representation of a documentational genesis*

The documentational genesis combines two interrelated processes: the instrumentalization (the teacher acting on resources), and the instrumentation process (resources supporting teacher's activity).

Taking into account teachers collective work leads us to articulate the documentational approach to the frame of the *community of practice* (Wenger 1998). A community of practice (CoP) is a human group presenting three main features: a *mutual engagement* of all its members, the active *participation* to the joint enterprise, and the *reification* of elements of practice (that is to say the production of things, results of the common practice, and recognized as a common wealth). This frame appears relevant for studying, instead of the interplay between a teacher and a set of resources (Fig. 1), the interplay between a group of teachers and sets of resources they are working on/with. We will call *community documentation* the process of reification carried out in teachers CoP.

## 3. OUR EXPERIMENTAL FIELD

Sésamath is a mathematics teachers association founded in 2001. Its kernel is constituted by about one hundred subscribers (mathematics teachers), sharing a set of principles inscribed in a charter[3]: common philosophy of public service, "mathematics for everybody"… Elected by this kernel, the Sésamath board regularly launches new projects groups for designing resources on a given theme, corresponding to teachers' special needs and interests (new subject in the curriculum, new textbook… see § 5.1, Fig. 3). These groups gather about five thousand teachers (a number of teachers belonging to several groups at the same time), working mainly at distance, via a platform and mailing lists. The group members can benefit from the assistance of employees and computer developers hired by the association (see fig. 3). All these groups present the features of CoP, of course not at the same level of development: the projects groups appear, at the beginning, as *potential* CoP (Wenger *et al.* 2002), while Sésamath kernel appears as a *maturing* one (ibidem). Sésamath, beyond its regular members, thus constitutes a *constellation of CoP* (ibidem), each of them sharing a same commitment for the principles of the association. A questionnaire (Sabra 2009), proposed in 2008 to 36 members of the Sésamath kernel during a training session[4], gave some answers (Table 1) to our first question (§1): why do teachers engage in Sésamath?

| What are, for your teaching, your sources of documentation? | Online resources (33/36); Resources that you have developed in previous years (35/36) |
|---|---|
| What are your professional reasons for engaging in Sésamath? | Training (16/36); Exchange of experience (23/36); Exchange of resources (13/36). |

**Table 1.** *Questionnaire to Sésamath, two questions (among 39) and their answers*

The first item confirms the Sésamath members' interest for online resources and evidences the place, for them, of resources *reuse* (revealing the importance, for Sésamath work, of exchanging, combining, modifying his/her own existing resources). The question item reveals that the main motivation, for joining Sésamath, is the *exchange of experience*, which evidences Sésamath role for its members' professional development.

Going further in our research required to follow the documentation work of one Sésamath project group on the field. We have chosen the group (named in the following DT10), created in 2009, aiming to design a digital textbook for the beginning (grade 10) of the French high school. This choice was motivated by two reasons: 1) after having designed textbooks for the French college (grade 7 to 9), Sésamath thus addresses, for the first time, the more complex mathematics of the grade 10 level; 2) after having designed textbooks quite classical (pdf, OpenOffice files that teachers could freely download and modify) with complement (spreadsheets and interactive applications), Sesamath aimed, with this project, to create a new type of *digital* textbook, that a teacher could appropriate and adapt to her own needs. We hypothesised that these new challenges (mathematical and technical) could enrich the documentation work of this group, and thus make it more interesting for our research.

Answering our initial questions needed also following the complete documentation work (both for Sésamath and for their own classes) of some DT10 members. For this article, we have chosen to present Anaïs' case. Anaïs is particularly engaged in DT10). She is 57 years old. After 15 years of various occupations (including some years as a member of a commune in the countryside), she came back to the university, achieved her mathematics studies, and got a position as a mathematics teacher. She has now 18 years of teaching experience (15 years for grade 10). Why does Anaïs engage in the DT10 project? How does her professional knowledge interact with individual and community documentation?

## 4. METHODOLOGY

Observing the individual and collective documentation is a complex task. It requires to take into account several conditions: long time observing to highlight the regularities; individual and collective observing; observing *in* and *out* of classroom; following both teachers' activities and resources. In some methodology of CoP observation, the researcher is engaged in practices (Jaworski 2009). In our case, we have just observed the practice as an outsider.

### 4.1 The reflexive investigation for observing individual documentation

For the observing of the Anaïs documentation, we adopted a methodology designed by Gueudet & Trouche (2010): the *reflexive investigation*. According to this methodology, the teacher is an essential actor in the collection of data. Among its methodological tools: interview at home, 'guided tour' (the teacher being the guide) and schematic representation of her *resource system* (SRRS); questionnaire about her vision of mathematics and of mathematics teaching; follow-up during several weeks including a logbook fulfilled by the teacher, collect of the resources, classroom observations. Anaïs became sick at the beginning of the follow-up, thus a direct classroom observation was impossible. We have adapted the methodology, in observing (and videotaping) Irvin,

Anaïs's substitute in class, who actually used Anaïs resources, and discussing afterwards with Anaïs on Irvin uses of her resources. In this paper, we only exploit data from the Anaïs' interview at home, from her comments on Irvin's video, from her SRRS and from a set of resources she designed for the teaching of functions (according to the discussion in DT10).

**4.2 Designing a methodology for observing community documentation**

Extending this methodology lying on *individual* reflexive investigation, new methodological tools have been designed for fostering teachers' reflection on their *collective* practice. Among them, an *agenda*, fulfilled by some members of the CoP, chosen for the role that they had in the project (for example, Adam was selected for its role identified in the mailing list as an effective coordinator of DT10). These agenda aims to identify and analyze, from different points of view of different actors, the effects of the *incidents* (something unexpected needing to reformulate the common goal or to reorganize the community documentation) occurring throughout the common project. We have also asked DT10 members a *schematic representation of collective interactions* (SRCI) in the case of Sésamath, and collect the mathematical resources that they designed.

The methodology takes also advantage of *natural data* which offers the experimental field. We thus exploit: the Anaïs online notebook; the DT10 mailing list, which offers discussions about project organisation, as well about mathematical, didactical and epistemological aspect of the community documentation work. During the period of our follow-up (11 months), 627 messages have been exchanged via the mailing list, involving 27 members (including Adam, Anaïs, John, Ben, Pierre and Henry, § 5). Anaïs has authored 111 messages among the 627. We particularly exploit in this paper the DT10 *thread of discussions* concerning the "mathematics functions", linked to an incident: a curriculum change occurring in the midst of DT10 work, and provoking the more intensive discussion.

**5. DATA ANALYSIS AND DISCUSSION**

We firstly analyze Sésamath and Anaïs' documentation, and then present the didactical interactions between Anaïs and DT10.

**5.1 Documentation work: DT10 and Anaïs**

Most of Anaïs's resources (courses, exercises, homework) are digital, stored in an external hard disc. This digital form facilitates the sharing with other teachers via USB key or Internet. The Anaïs' resource system is strongly articulated with Sésamath resources (see Anaïs SRRS, Fig. 2): emails, students' sheets and other resources exchanged with Sésamath seem to have a major role in her documentation. This *osmosis* between Anaïs and Sésamath resources facilitates Anaïs' participation in DT10 work.

Anaïs's documentation is conditioned by three factors: "curriculum changing, institutional recommendations and classroom general level" (Anaïs' interview). For example, her teaching of "function" is strongly related to curriculum change. She said that she lived "*two different spirits*" of teaching functions: "it was, before, more guided and now the curriculum recommends open questions... so that students have more initiative" (Anaïs' interview). She seems very sensitive to the students' level: this year, I had a low level class. I will modify [the documents of the previous year] by adapting them to the level of my class. But maybe, the next year I will have good students, I will reuse the document in the present form" (Anaïs' interview).

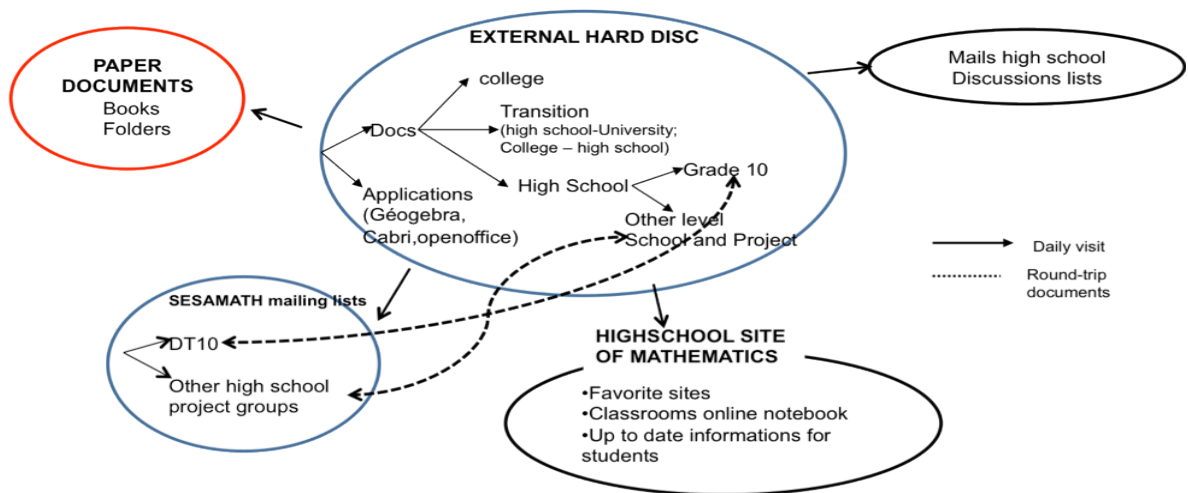

**Figure 2.** *Anaïs' SRRS, our translation and schema as close as possible of the original*

DT10 documentation is also sensitive to curriculum changes and keeps Sésamath resources as a general background support (particularly when an incident occurs), as it appears in the following extract of the DT10 mailing list:

> **John:** The curriculum has profoundly changed and we must rebuild the chaptering. This chaptering is very important because it is common to several Sésamath projects that involve the grade 10. Thus, I propose to do a general RESET.
>
> **Henry:** We will create a general mailing list for all Sésamath projects linked to the secondary level as there issues intersect. We will thus share resources and general links that might be useful to all.
>
> **John:** This list must first allow everyone to have a clear idea of what is done in different projects.

DT10 and Anaïs's documentation do not answer to the same constraints: DT10 has to move forward systematically with other high school projects, while Anaïs is attracted by the level of her pupils, particularly in the new curriculum spirit. But both documentations are conditioned by the institutional recommendations, which constitute the joint constraint for designing the resources.

Anaïs's interactions with DT10 are linked to her general interest for collective matters, and particularly in this project. Anaïs underlines, in her SRCI for Sésamath (Fig. 3), that each project is constituted by teachers which are interested the project's topic. More precisely, her participation resulted from a professional interest: participating in DT10, constitutes, for her, an occasion for searching ideas for her own classroom course. She confirms this idea in the interview, illustrating the productive aspects of instrumentalisation process (§ 2): "for example, for the chapter "Variations and Extrema of functions" [named in the following 2N2], I don't have constructed the lesson yet. In fact, I will transform it [the lesson designed in DT10] and I will reuse it for my classroom. I will adapt it for my classroom".

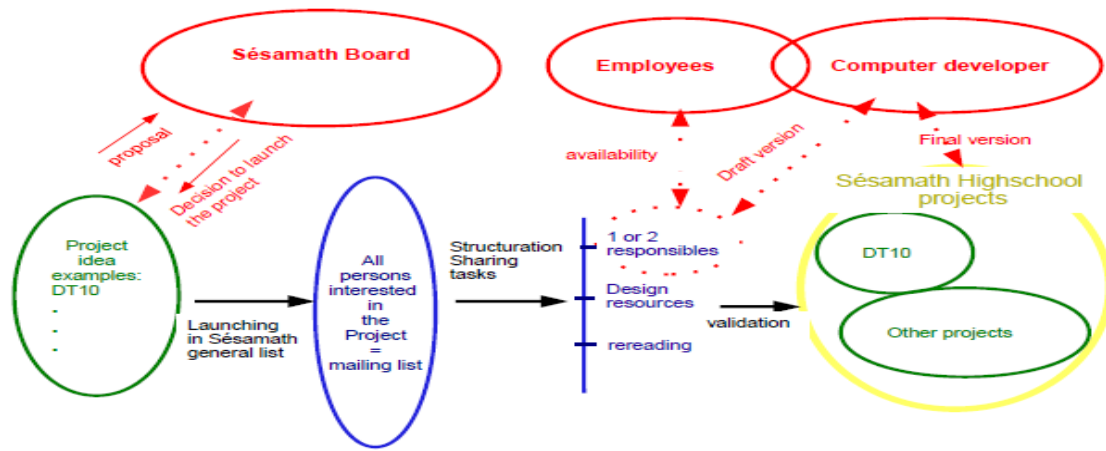

**Figure 3.** *Anaïs's SRCI, our translation (see comments on project groups §3)*

We analyze in the following section the didactical aspects of the interactions between Anaïs and DT10 members.

## 5.2 The didactical interaction between Anaïs and DT10 members

The didactical interactions between Anaïs and other DT10 members take place, mainly, through the mailing list of the project. Most of the didactical discussion, in the thread "Create textbook chapters", concerns the theme "mathematics functions". This discussion came under the PCK category "instructional strategies for teaching the functions", as it appears in the following extract of DT10 mailing list about the 2N2 chapter (evidencing Adam's role as coordinator, see § 4.2):

**Ben:** We need to cut up 2N2 chapter. What are your proposals?
**Pierre:** … How far will we go with the concept of variations? My proposal: 1) the notion of functions variation (from a curve); 2) Maximum and minimum of a function, 3) graph and table of variation; 4) comparison of numbers.
**Anaïs:** My proposals: 1) functions variations and graphical reading... 2) variations and calculations: square function ...inverse and linear functions, 3) Maximum and minimum of a function, 4) table of variation.
**Adam:** I have a slight preference for Pierre's proposal ... because I think that the table of variation is simpler to understand than computations for comparing numbers ... I think that slowly increasing the difficulty is a good strategy.

Anaïs's chapter happens following another structure. She first introduces the concept from what she calls a "start-up activity", with the aim of presenting the semantic values of the "mathematics functions" terminology; then, "the course" where she presents the main definitions (and related examples); at the end, the "series of exercises", and "homework and assessment". This structure appears in her online notebook (table 2).

| Monday 21/09/2009 | Constructing, from a sheet of paper, a box whose volume is the bigger one |
|---|---|
| Friday 25/09/2009 | Definitions: set of a function definition, variable, image, and antecedent. Table of values and graph using a calculator. |
| Monday 28/09/2009 | Area of a rectangle with a given perimeter. Table of values and graph of *reference* functions (that are inverse function, square function, square root function…) |
| Thursday October 1st | Different ways of saying that f (a) = b |
| Friday 2/10/2009 | Exercises on image of an interval; first assessment. |

**Table 2.** *An extract of the Anaïs' online notebook*

In the thread of discussion "Extremum 2N2", that concerns designing an exercise for 2N2, Anaïs initiates didactical discussions about the terminology:

**Anaïs:** f(x) is smaller than (or the image f(x) is smaller than?). It seems...

**John:** For M being a maximum, it must be both an upper bound and an image by f. That is to say, for all x belonging to a given interval I, f(x) ≤4.5, then 4.5 is an upper bound. And it is a maximum if it is also an image of some x belonging to I.

**Adam:** I propose this formulation: for every real x belonging to I, f(x) is smaller than, or equal to, f(3), which equals 4.5.

**Anaïs:** Either like this: for every x in the interval I, f(x) ≤4.5 with 4.5=f(3).

Also, when she discovered Irvin's video, she has particularly commented the language used by the teacher:

**Anaïs:** the expressions ... he says "the straight line 3x +1", or the straight line Y. What does that mean? He should say the line whose equation is y = 3x+1.

**Me:** why does it bother you? How this is a problem?

**Anaïs:** I think we should be highly accurate, even when we talk. Being vague or loosely makes fuzzy pupils' mind, especially when introducing new concepts...

We remark a difference between DT10 (more precisely DT10 members taking in charge this question) and Anaïs in the instructional strategies for teaching function. DT10 constructs a "function" concept from very simple tools (reading of curves and table of values) then it moves to more complex tools like calculation and articulation between different types of representation (algebraic, graphical). Anaïs takes care, in her documentation, of the terminology of mathematical concept and its use. This *gap* is an opportunity for Anaïs to discuss her view with other DT10 members.

On another subject, the design of a problem-situation reveals a *convergence* between Anaïs and other DT10 members documentation: Anaïs fosters the place of conjecture and experimentation. She presents thus a resource (interview):

> **Anaïs:** I like to ask them first to experiment with a calculator... try first with the calculator, guess the number of solutions ... for example f (x) = 0 they speculate with the calculator and ... Once they guess, after they do the proof by calculation. **But I think it is important to first think about from free explorations.**

The thread of discussion "test_2N2" reveals, within DT10, a shared point of view on problem-situation. The problems designed consist in modelling geometrical situations with a function, graphing this function with a calculator and elaborating conjectures before computing and searching for evidence.

## 6. CONCLUSION

We have presented in this paper a theoretical approach and a methodology aiming to analyze the didactical interactions at stake within a *community documentation work*. We have carried out this approach and this methodology in the case of a documentation project, DT10, of an online teachers association, focusing on Anaïs, a member of this association.

Anaïs and DT10 have a common interest in collaborating for designing resources and collectively facing *incidents* (like changes in the curriculum). Anaïs has also an individual interest linked to her documentation needs: discussing the mathematics terminology and language, and its values in the documentation work. When Anaïs perceives that didactical discussions are part of her interest, she participates strongly in the community documentation.

Three factors seem to stimulate the active participation of a teacher in the community documentation: osmosis between her resource system and the CoP resources; a gap between the teacher and the CoP in the strategies for teaching a subject; a shared interest in the subject of discussion that is the origin of gap. Following an episode of Anaïs's active participation in the community documentation, we identified an aspect of her professional knowledge: the introduction of a concept has to be based on a precise introduction of terminology and language associated, and its uses.

This study revealed a question that deserves more deepening: how identifying the effects of the community documentation on the teacher professional knowledge in the case of convergence between the individual and community documentation?

## REFERENCES

Grossman, P. (1990), *The Making of a Teacher*, New York: Teachers College Press.

---

[1] http://www.Sesamath.net/

[2] The choice of vocabulary intended to match the terminology of document management research. According to Pédauque (2006), "A document is not anything, but anything can become a document, as soon as it supplies information, evidence, in short, as soon as it is authoritative." (p. 12, our translation).

[3] See: http://www.Sésamath.net/association.php?page=asso_charte.

[4] This training session was organized by both Sésamath and the National Institute for Pedagogical Research.